# A theorem on circle configurations

*Jerzy Kocik*

jkocik@siu.edu
Mathematics Department, SIU-C, Carbondale, IL

**Abstract:** A formula for the radii and positions of four circles in the plane for an arbitrary linearly independent circle configuration is found. Among special cases is the recent extended Descartes Theorem on the Descartes configuration and an analytic solution to the Apollonian problem. The general theorem for *n*-spheres is also considered.

## 1. Introduction

There is a family of problems of various degrees of generality concerning circles and, in general, $(n-1)$-dimensional spheres in $\mathbf{R}^n$. Among the most famous is the problem of Apollonius of Perga: find a circle simultaneously tangent to three given circles (Fig. 1.1a); and the problem of Descartes: given three pair-wise tangent circles, find a fourth tangent circle to each of them. A solution to the latter problem — four mutually tangent circles — is known as a Descartes configuration (Fig. 1.1b). The associated Descartes formula relates the sizes of the circles in a peculiarly elegant way:

$$(a + b + c + d)^2 = 2(a^2 + b^2 + c^2 + d^2), \tag{1.1}$$

where $a=1/r_1$, $b=1/r_2$, etc. are the curvatures of the corresponding circles (*reciprocals* of radii). For the history of this formula, rediscovered a number of times, and its higher-dimensional generalizations, see [Des, Ped3, Coxt69, Sodd, Gos].

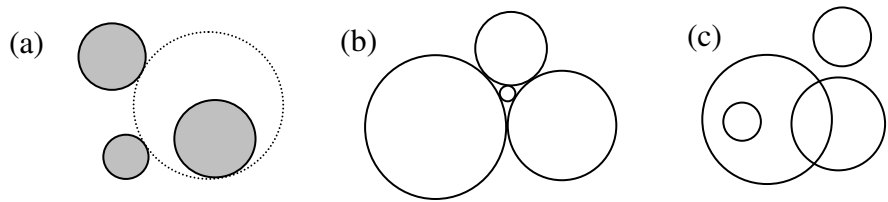

**Figure 1.1**: Four circles ; (a) as a solution to an Apollonian problem;
(b) in Descartes configuration (c) in a general configuration

It's easy to observe that the quadratic formula 1.1 may be written in a matrix form:

$$\begin{bmatrix} b_1 & b_2 & b_3 & b_4 \end{bmatrix} \begin{bmatrix} -1 & 1 & 1 & 1 \\ 1 & -1 & 1 & 1 \\ 1 & 1 & -1 & 1 \\ 1 & 1 & 1 & -1 \end{bmatrix} \begin{bmatrix} b_1 \\ b_2 \\ b_3 \\ b_4 \end{bmatrix} = 0, \tag{1.2}$$

i.e., $\mathbf{b}^T D \mathbf{b} = 0$, where $\mathbf{b}$ is the vector of four curvatures and where $D$ has been termed the "Descartes quadratic form".

The Darboux formula has recently been extended to matrix form so that it includes the position of the centers of the circles [LMW]. The proof goes via hyperbolic geometry.

In this paper we go beyond this "tangential" canon and present a theorem on an (almost) arbitrary configuration of four circles. It provides a formula relating their radii (curvatures) and positions, given some mutual relations, like tangency, orthogonality, distances, etc. A generalization to $n$-dimensional spheres is also given.

Our proof relies on the fact that circles ($n$-spheres) may be mapped to the vectors of a Minkowski space, a discovery made – to our knowledge – by D. Pedoe. We use the induced inner product in the dual space.

The theorem presented here provides a simple but powerful tool for generating "Descartes-like theorems" for different types of configurations, the Descartes case being just one example, as well as simply to solve particular configurations.

## 2. Geometry of Circles — the Pedoe map

Circles in the plane may be viewed as vectors of a Minkowski space with a pseudo-Euclidean inner product. This beautiful fact ties the ancient geometry of circles and the modern geometry of space-time. Let us recall the development of this idea.

**A. Remark on the history of the Pedoe product:** In 1826, Jakob Steiner defined the *power* of a point $P$ with respect to a circle $C$:

$$P * C = d^2 - r^2, \qquad (2.1)$$

where $r$ = radius, $d$ = distance from P to the center of the circle. Consult Fig 2.2a for the motivation: the product of segment lengths $PA \cdot PB$ does not depend on the choice of the line through P. Choosing the line tangent to the circle gives expression (2.1). In 1866, G. Darboux [Dar] generalized this notion to the *power* of a pair of circles, which we will call a ***Darboux product***:

$$C_1 * C_2 = d^2 - r_1^2 - r_2^2. \qquad (2.2)$$

If the circles intersect, the product equals $C_1 * C_2 = r_1 r_2 \cos\varphi$, where $\varphi$ is the angle made by the circles (see Figure 2.2a). Let us note that in the case of more distant circles, the product $C_1 * C_2$ equals the square of the segment constructed in Figure 2.2c.

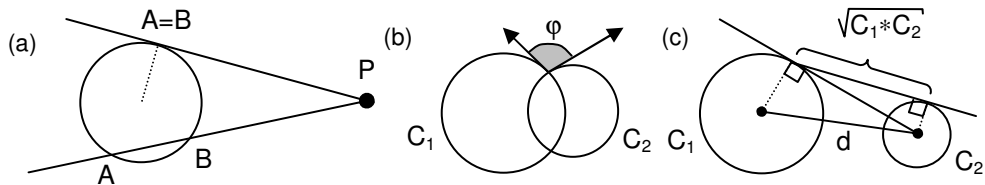

**Figure 2.2:** (a) Power of a point. Geometric interpretation of the Darboux product of (a) intersecting circles; (c) distant circles.

The equation of the circle of radius $r$ centered at $(x_o, y_o)$ has the form:

$$x^2 + y^2 - 2xx_o - 2yy_o + c = 0, \qquad (2.3)$$

where $c = x_o^2 + y_o^2 - r^2$. H. Cox wrote the Darboux product of two circles with centers $(x_1, y_1)$ and $(x_2, y_2)$ and radii $r_1$ and $r_2$, respectively, in terms of the coefficients of the corresponding equations:

$$C_1 * C_2 = c_1 + c_2 - 2x_1x_2 - 2y_1y_2 \tag{2.4}$$

([Cox]). The crucial observation was made in 1970 by D. Pedoe, who realized that (2.4) may be interpreted as an inner product [Ped2]. Although (2.4) contains linear terms, $c_1$ and $c_2$, he pointed out that the equation of a circle is invariant with respect to multiplication by a scalar, and its general form is

$$a(x^2 + y^2) - 2px - 2qy + c = 0. \tag{2.5}$$

One can thus introduce a scalar product $\langle C_1, C_2 \rangle$ as

$$2 \langle C_1, C_2 \rangle = a_1 c_2 + c_1 a_2 - 2p_1 p_2 - 2q_1 q_2. \tag{2.6}$$

Note that when $a_1$ and $a_2$ are chosen to be $a_1 = a_2 = 1$, the scalar product coincides with (half) the Darboux product (2.4). (The factor of 2 is introduced for later convenience).

**B. Pedoe map.** The purpose of the previous subsection was to establish that the discovery of the pseudo-Euclidean geometry of circles should be credited to D. Pedoe, even though he did not follow the path of Minkowski geometry. Let us formalize this discovery and introduce some terminology.

**Definition 2.1:** By a **standard isotropic** 4-dimensional Minkowski space we understand a linear real space $M \cong \mathbf{R}^4$ with an inner product $\langle , \rangle$ given by the *metric matrix*

$$g = \begin{bmatrix} 0 & \frac{1}{2} & 0 & 0 \\ \frac{1}{2} & 0 & 0 & 0 \\ 0 & 0 & -1 & 0 \\ 0 & 0 & 0 & -1 \end{bmatrix} \tag{2.7}$$

For vectors $\mathbf{v} = [v_1, v_2, v_3, v_4]$ and $\mathbf{w} = [w_1, w_2, w_3, w_4]$, we have

$$\langle \mathbf{v}, \mathbf{w} \rangle = \tfrac{1}{2} v_1 w_2 + \tfrac{1}{2} v_2 w_1 - v_3 w_3 - v_4 w_4 \qquad \textit{[inner product]}$$

$$|\mathbf{v}|^2 = v_1 v_2 - v_3^2 - v_4^2. \qquad \textit{[norm squared]}$$

The basis in which the above inner product is expressed will be called the **standard isotropic basis**.

**Definition 2.2: The Pedoe projective map** is a map from the set $\Omega$ of circles in the plane to the rays in the standard isotropic Minkowski space **M** (projective space *PM*)

$$\pi: \Omega \to PM$$

so that the circle $C_r(x_0, y_0)$ with center at $(x_0, y_0) \in \mathbf{R}^2$ and radius $r \in \mathbf{R}$ is mapped to a ray $\pi(C)$ in *M* spanned by a vector:

$$\pi(C) = \text{span} \{ [ 1, x_0^2 + y_0^2 - r^2, x_0, y_0 ]^T \}.$$

and a line *L* (circle of null curvature) given by equation $ax + by = c$ is mapped to

$$\pi(L) = \text{span} \{ [ 0, c/2, a, b ]^T \}.$$

Note that the Pedoe map may be extended to points (**improper circles**) so that point $P = (x, y)$ has image $\pi(P)$ = span { $[ 1, x^2+y^2, x, y]^T$ }, which is a ray in the light cone, as $|\pi(v)|^2 = 0$, $\forall v \in \pi(P)$. It is easy to see that the rays that represent proper circles are space-like, i.e., they lie in *M* outside the light cone. For them, we may restrict the Pedoe map to the unit vectors.

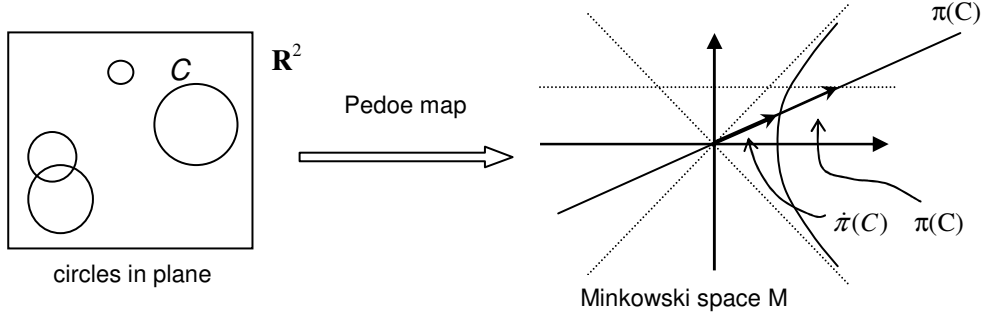

**Figure 2.2:** The Pedoe map carries circles to rays in the Minkowski space

**Definition 2:** The **Pedoe special map** is a specification of $\pi$: $\Omega \to$ PM, i.e., is map

$$\dot{\pi}: \Omega \to M$$

that carries circles to unit space-like vectors of *M* so that span $\dot{\pi}(C) = \pi(C)$, $|\dot{\pi}(C)|^2 = -1$. To remove ambiguity, we require that the first component of $\dot{\pi}(C)$ be nonnegative. Vector $\dot{\pi}(C)$ will be called the **Pedoe vector** of circle *C*. In particular, a circle of radius $r \neq 0$ centered at $(x,y)$ is represented by the Pedoe vector

$$\dot{\pi}(C) = \begin{bmatrix} b \\ \bar{b} \\ \dot{x} \\ \dot{y} \end{bmatrix} \equiv \begin{bmatrix} 1/r \\ (x^2 + y^2 - r^2)/r \\ x/r \\ y/r \end{bmatrix}. \tag{2.8}$$

Its entries will be called:

$$\begin{aligned}
\text{circle curvature:} \quad & b = 1/r \\
\text{circle co-curvature:} \quad & \bar{b} = (x_0^2 + y_0^2 - r^2)/r \\
\text{reduced position:} \quad & \dot{x} = x/r, \quad \dot{y} = y/r
\end{aligned} \tag{2.9}$$

(Note that the second term of the Pedoe vector is determined by the other three and by the requirement of normalization).

**Remark on the geometric meaning:** The *curvature b* has the standard meaning. The *co-curvature* $\bar{b}$ happens to equal the curvature of the circle that is the image of *C* via inversion in the unit circle centered at the origin. As for the reduced position ($\dot{x}$, $\dot{y}$), one can think of it as the position of an "effective" circle one would regard as the circle viewed from the origin, in the belief that it has radius equal to 1.

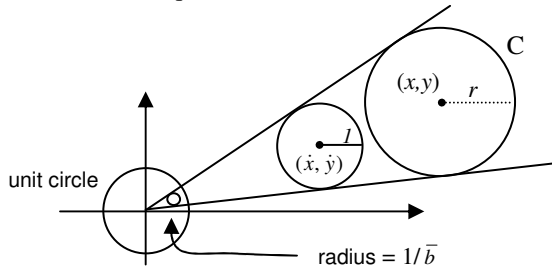

**Figure 2.3:** The geometric meaning of reduced position ($\dot{x}$, $\dot{y}$) and of co-radius $\bar{b}$

**Remark:** Thus, the following vectors represent the same circle:

$$C \to \begin{bmatrix} a \\ c \\ p \\ q \end{bmatrix} \sim \begin{bmatrix} 1 \\ x_0^2 + y_0^2 - r^2 \\ x_0 \\ y_0 \end{bmatrix} \sim \dot{\pi}(C) = \begin{bmatrix} b \\ \bar{b} \\ \dot{x}_0 \\ \dot{y}_0 \end{bmatrix} = \begin{bmatrix} 1/r \\ (x_0^2 + y_0^2 - r^2)/r \\ x_0/r \\ y_0/r \end{bmatrix},$$

where the first vector is the set of coefficients of Eq. (2.6), the second is scaled by setting $a=1$, and the third, scaled by the radius of the circle, is the Pedoe vector. The first parameterization of vectors has the advantage that choosing $a = 0$ admits lines into the formalism as special circles. Points $(x,y)$ go to $[1, x^2+y^2, x, y]$, corresponding to light-like rays; hence they do not permit Pedoe vectors.

**Theorem 2.3 [Pedoe]:** The special Pedoe map is an injection into the unit hyperboloid in $M$ and the corresponding pseudo-Euclidean inner product is related to the Darboux product of circles:

$$\langle \pi(C_1), \pi(C_2) \rangle = \tfrac{1}{2} C_1 * C_2 / r_1 r_2 \qquad (2.10)$$

The above expression will be called the **Pedoe inner product** of circles.

**Proposition 2.4:** Let $\mathbf{C}_i = \dot{\pi}(C_i)$, $i=1,2$, denote unit vectors representing a pair of circles. Then:

a) $|\mathbf{C}_i|^2 = -1$ (space-like, unit vector) for any circle

b) $\langle \mathbf{C}_1, \mathbf{C}_2 \rangle = \cos\varphi$, if $C_1 \cap C_2 \neq \emptyset$

c) $\langle \mathbf{C}_1, \mathbf{C}_2 \rangle = +1$ if $C_1$ and $C_2$ are tangent externally 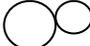

d) $\langle \mathbf{C}_1, \mathbf{C}_2 \rangle = -1$ if $C_1$ and $C_2$ are tangent internally 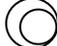

e) $\langle \mathbf{C}_1, \mathbf{C}_2 \rangle = 0$ if $C_1$ and $C_2$ are mutually orthogonal 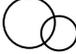

**Proof:** Simple verification.

**Remark on basis: isotropic versus orthonormal.** Until now we have dealt with the isotropic basis and non-diagonal metric tensor $g$. But one may want to choose an orthonormal basis in which the metric tensor is diagonalized,

$$g = \begin{bmatrix} 0 & 1/2 & 0 & 0 \\ 1/2 & 0 & 0 & 0 \\ 0 & 0 & -1 & 0 \\ 0 & 0 & 0 & -1 \end{bmatrix} \sim g_0 = \begin{bmatrix} 1 & 0 & 0 & 0 \\ 0 & -1 & 0 & 0 \\ 0 & 0 & -1 & 0 \\ 0 & 0 & 0 & -1 \end{bmatrix},$$

so that the the Pedoe inner product signature $\text{sign}(g) = (+,-,-,-)$ is more conspicuous. The scalar product is here simply $\langle \mathbf{v}, \mathbf{w} \rangle = v_1 w_1 - v_2 w_2 - v_3 w_3 - v_4 w_4$. The Pedoe vector of a circle $C_r(x,y)$ is now:

$$\dot{\pi}(C) = \begin{bmatrix} \frac{1+\rho^2 - r^2}{2r} \\ \frac{1-\rho^2 + r^2}{2r} \\ \dot{x} \\ \dot{y} \end{bmatrix} = \begin{bmatrix} \frac{b^2 + \dot{\rho}^2 - 1}{2b} \\ \frac{b^2 - \dot{\rho}^2 + 1}{2b} \\ \dot{x} \\ \dot{y} \end{bmatrix}$$

where $\rho^2 = x_0^2 + y_0^2$.

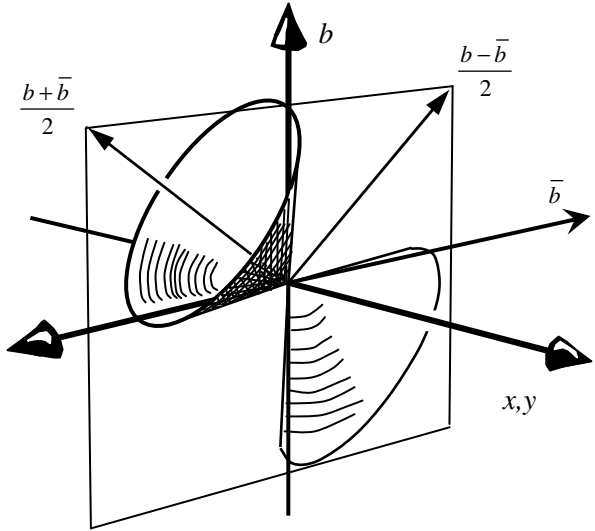

**Figure 2.4:** Isotropic basis versus orthonormal basis in the Minkowski space for circles

Although the orthonormal basis side is more "suggestive" for those accustomed to relativistic physics (see e.g. [Sod]), the path to this form leads from the Darboux product in a way that is not as intuitive and clear with respect to the isotropic basis and coordinates. Figure 2.4 shows the bases and corresponding variables symbolically, squeezing 4 dimensions into a 3-dimensional picture. One of the axes carries variables $x$ and $y$ (the position of the circle's center). (The figure represents exactly the special case of n=1, where spheres are pairs of points on a line).

$$\begin{bmatrix} 1 \\ x_0^2 - r^2 \\ \boldsymbol{x}_0 \end{bmatrix} \sim \begin{bmatrix} 1/r \\ (x_0^2 - r^2)/r \\ \boldsymbol{x}_0/r \end{bmatrix} = \begin{bmatrix} b \\ \bar{b} \\ \dot{\boldsymbol{x}}_0 \end{bmatrix} \quad \underset{\text{basis}}{\overset{\text{change of}}{\Longleftrightarrow}} \quad \begin{bmatrix} (b+\bar{b})/2 \\ (b-\bar{b})/2 \\ \dot{\boldsymbol{x}}_0 \end{bmatrix} = \begin{bmatrix} \dfrac{1+x_0^2-r^2}{2r} \\ \dfrac{1-x_0^2+r^2}{2r} \\ \dfrac{\boldsymbol{x}_0}{r} \end{bmatrix} \sim \begin{bmatrix} 1+x_0^2-r^2 \\ 1-x_0^2-r^2 \\ 2\boldsymbol{x}_0 \end{bmatrix}$$

norm = −1          norm = −1

$$\text{here } g = \begin{bmatrix} 0 & \tfrac{1}{2} & 0 \\ \tfrac{1}{2} & 0 & 0 \\ 0 & 0 & -\mathbf{1} \end{bmatrix} \qquad \text{here } g' = \begin{bmatrix} 1 & 0 & 0 \\ 0 & -1 & 0 \\ 0 & 0 & -\mathbf{1} \end{bmatrix}$$

**Figure 2.5:** Change of basis

Circles correspond to the rays (vectors) outside the cone. Points in $\mathbf{R}^n$, which may be viewed as circles of null radius, correspond to light-like rays. Similarly, lines, which may be viewed as circles of zero curvature, are rays in the horizontal hyper-plane $b = 0$.

# 3. The theorem on Circle configuration

Consider four circles $C_1, \ldots, C_4$ represented by Pedoe unit vectors (see (2.10)). Define for this set of circles a *configuration matrix f* as the Grammian of the vectors $C_i$, that is:

$$f_{ij} = \langle \mathbf{C}_i, \mathbf{C}_j \rangle, \tag{3.1}$$

where $\langle \cdot, \cdot \rangle$ is the Pedoe inner product. On the other hand, we may build a "*data matrix*" as the collective representation of the circles (columns):

$$A = [\mathbf{C}_1 \mid \mathbf{C}_2 \mid \mathbf{C}_3 \mid \mathbf{C}_4]. \tag{3.2}$$

Now we are ready to state and prove our theorem in a shockingly simple way.

**Theorem 3.1 (Circle Configuration Theorem):** If circles $C_1, \ldots, C_4$ are linearly independent, then

$$\boxed{AFA^T = G} \tag{3.3}$$

where the matrix $G$ is the inverse of the matrix of Minkowski metric, $G = g^{-1}$, and F is the inverse of the configuration matrix, $F = f^{-1}$.

**Proof:** Definition (3.1) of the configuration matrix $f$ may be written in cumulative matrix form as

$$f = A^T g A. \tag{3.4}$$

Since $A$ is invertible due to the linear independence of the four circles, we may take the inverse of both sides

$$f^{-1} = A^{-1} g^{-1} (A^T)^{-1}.$$

Now, by multiplying both sides on the left and right by $A$ and $A^T$, respectively, we get (3.3) as desired.

Let us try to understand the benefits of this equation. First, we introduce four vectors that represent the data on the circles in a way *dual* to what we used so far, namely we consider

$$\boldsymbol{b} = \begin{bmatrix} b_1 \\ b_2 \\ b_3 \\ b_4 \end{bmatrix} \quad \dot{\boldsymbol{x}} = \begin{bmatrix} \dot{x}_1 \\ \dot{x}_2 \\ \dot{x}_3 \\ \dot{x}_4 \end{bmatrix} \quad \dot{\boldsymbol{y}} = \begin{bmatrix} \dot{y}_1 \\ \dot{y}_2 \\ \dot{y}_3 \\ \dot{y}_4 \end{bmatrix} \quad \bar{\boldsymbol{b}} = \begin{bmatrix} \bar{b}_1 \\ \bar{b}_2 \\ \bar{b}_3 \\ \bar{b}_4 \end{bmatrix}.$$

Vector $\boldsymbol{b}$ will be called the *curvature vector* of the circle configuration, $\bar{\boldsymbol{b}}$ the *co-curvature vector*, and $\dot{\boldsymbol{x}}$, and $\dot{\boldsymbol{y}}$ are the *reduced position vectors*. The data matrix constructed by superposing these columns is $A^T$

$$A^T = [\boldsymbol{b} \mid \bar{\boldsymbol{b}} \mid \dot{\boldsymbol{x}} \mid \dot{\boldsymbol{y}}] = \begin{bmatrix} b_1 & \bar{b}_1 & \dot{x}_1 & \dot{y}_1 \\ b_2 & \bar{b}_2 & \dot{x}_2 & \dot{y}_2 \\ b_3 & \bar{b}_3 & \dot{x}_3 & \dot{y}_3 \\ b_4 & \bar{b}_4 & \dot{x}_4 & \dot{y}_4 \end{bmatrix}.$$

We shall also need the explicit form of the inverse of the matrix of Minkowski metric:

$$g = \frac{1}{2}\begin{bmatrix} 0 & 1 & 0 & 0 \\ 1 & 0 & 0 & 0 \\ 0 & 0 & -2 & 0 \\ 0 & 0 & 0 & -2 \end{bmatrix} \quad \Rightarrow \quad G \equiv g^{-1} = \begin{bmatrix} 0 & 2 & 0 & 0 \\ 2 & 0 & 0 & 0 \\ 0 & 0 & -1 & 0 \\ 0 & 0 & 0 & -1 \end{bmatrix}. \quad (3.5)$$

Now, we see that Theorem 3.1 may equivalently be stated:

**Corollary 3.2:** Let $\mathbf{v}_i$, $i = 1,\ldots 4$, be one of the four vectors: $\mathbf{b}, \overline{\mathbf{b}}, \dot{\mathbf{x}}$, or $\dot{\mathbf{y}}$. Then

$$\mathbf{v}_i^T F \mathbf{v}_j = G_{ij}$$

In particular, for the vector of curvatures and reduced positions we have these handy quadratic formulas (given also in indexed form):

$$\begin{aligned} \dot{\mathbf{x}}^T F \dot{\mathbf{x}} &= -1 \quad \text{or} \quad \dot{x}_i F_{ij} \dot{x}_j = -1, \\ \dot{\mathbf{y}}^T F \dot{\mathbf{y}} &= -1 \quad \text{or} \quad \dot{y}_i F_{ij} \dot{y}_j = -1, \\ \mathbf{b}^T F \mathbf{b} &= 0 \quad \text{or} \quad b_i F_{ij} b_j = 0, \end{aligned} \quad (3.6)$$

where the indices label the four circles, $i,j = 1,\ldots 4$, and where summation over repeated indices is understood. The last equation is a Descartes-like formula —generalized to arbitrary independent circle configurations.

Let us now look at the theorem in action.

**Example 3.3 [the circle of inversive symmetry]:** Suppose there are given three pairwise externally tangent circles of curvatures $b_1=a$, $b_2=b$, and $b_3=c$. A fourth circle of curvature $b_4=d$ is orthogonal to each of the three.

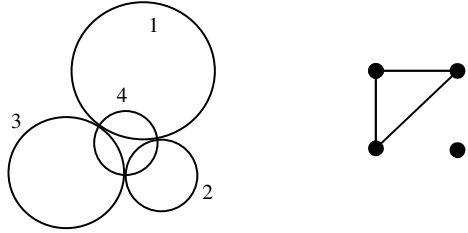

**Figure 3.1:** Circle configuration for Example 3.3. The symbol to the right of the configuration will be explained later.

We set the configuration matrix $F$ and calculate its inverse:

$$f = \begin{bmatrix} -1 & 1 & 1 & 0 \\ 1 & -1 & 0 & 0 \\ 1 & 1 & -1 & 0 \\ 0 & 0 & 0 & -1 \end{bmatrix}, \quad F \equiv f^{-1} = \frac{1}{2}\begin{bmatrix} 0 & 1 & 1 & 0 \\ 1 & 0 & 1 & 0 \\ 1 & 1 & 0 & 0 \\ 0 & 0 & 0 & -2 \end{bmatrix}.$$

The circle configuration theorem gives here

$$\begin{bmatrix} b_1 & b_2 & b_3 & b_4 \\ \overline{b}_1 & \overline{b}_2 & \overline{b}_3 & \overline{b}_4 \\ \dot{x}_1 & \dot{x}_2 & \dot{x}_3 & \dot{x}_4 \\ \dot{y}_1 & \dot{y}_2 & \dot{y}_3 & \dot{y}_4 \end{bmatrix} \begin{bmatrix} 0 & 1 & 1 & 0 \\ 1 & 0 & 1 & 0 \\ 1 & 1 & 0 & 0 \\ 0 & 0 & 0 & -2 \end{bmatrix} \begin{bmatrix} b_1 & \overline{b}_1 & \dot{x}_1 & \dot{y}_1 \\ b_2 & \overline{b}_2 & \dot{x}_2 & \dot{y}_2 \\ b_3 & \overline{b}_3 & \dot{x}_3 & \dot{y}_3 \\ b_4 & \overline{b}_4 & \dot{x}_4 & \dot{y}_4 \end{bmatrix} = 2\begin{bmatrix} 0 & 4 & 0 & 0 \\ 4 & 0 & 0 & 0 \\ 0 & 0 & -2 & 0 \\ 0 & 0 & 0 & -2 \end{bmatrix}.$$

(We multiplied both sides by the factor 2). The formula for curvatures (an analog of the Descartes formula) may be read off from the top left entries:

$$d^2 - ab - ac - ca = 0,$$

which resolves into the well-known formula:

$$d = \pm\sqrt{ab + bc + ca}\,,$$

or $r_4^2 = r_1 r_2 r_3 /(r_1+r_2+r_3)$. The theorem predicts also the position of the center of the orthogonal circle in terms of the data for the other circles:

$$\dot{x}_4^2 = \dot{x}_1\dot{x}_2 + \dot{x}_2\dot{x}_3 + \dot{x}_3\dot{x}_1 - 1 \quad \text{or} \quad x_4^2 = \frac{r_1 x_2 x_3 + r_2 x_3 x_1 + r_3 x_1 x_2 - r_1 r_2 r_3}{r_1 + r_2 + r_3},$$

and similarly for $y_4$. Thus the equations fully determine the fourth circle.

An immediate implication of the theorem is this geometric fact:

**Corollary 3.4:** There is no configuration of four mutually perpendicular circles.

**Proof:** The configuration matrix for such an arrangement of circles is $F = I = \text{diag}(1,\ldots,1)$. The "master equation" (3.3) becomes thus

$$AA^\text{T} = G. \tag{3.6}$$

This is impossible, for the left side is positive definite, while the right side is not. Restating it: Equation (3.3) is a congruence relation, but $F=I$ is not congruent to $G$ since they have different Sylvester's moments, $(+, +, +, +)$ versus $(-, +, +, +)$.

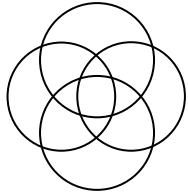

**Figure 3.2:** A failed attempt to draw four mutually orthogonal circles.

This of course is a special case of a more general statement:

**Corollary 3.5:** If the configuration matrix of a hypothetical arrangement of four circles has eigenvalues different from one strictly positive and three strictly negative, then the configuration cannot be realized.

**Remark on the interpretation of the formulas:** The formula (3.4) has a clear geometric meaning: The four circles $C_1,\ldots,C_4$ — if linearly independent — define a new basis in Minkowski space. Then $A^\text{T} g A = f$ represents the congruence of matrices $g \sim f$, i.e., the transformation of the Minkowski metric matrix under a change of basis, where

- $g$ = matrix of the Pedoe quadratic form in the standard (Pedoe) basis,
- $f$ = matrix of the same Pedoe quadratic form, expressed in the basis made by the four circles.

Now, in order to organize this equation in terms of the curvature vector **b** and (reduced) positions, which correspond to the rows of $A$, we need to move to the *dual* Minkowski space. The "master formula" of Theorem 3.1 does exactly this: it represents congruency of the induced metric matrices $F$ and $G$ in the *dual space*. The columns of matrix $A$ originate as vectors, while the rows, including **b** = $[b_1, b_2, b_3, b_4]$, correspond to covectors, representing the basis dual to that of the basis given by the circles.

# 4. More examples and new configuration families

The Descartes configuration is but one of many families of configurations. To exemplify this point, let us have a look at two such families (they extend the first two examples of the last section).

**Example 1: Descartes configuration revisited.**

Among the special cases of Theorem 3.1 is the original Descartes formula and its extension, discovered by J. Lagarias *et al*. [LMW]. Let us look at some details. Four circles are said to be in a *Descartes configuration* if all pairs of circles are mutually tangent at distinct points (see Figure 4.1 for possible arrangements).

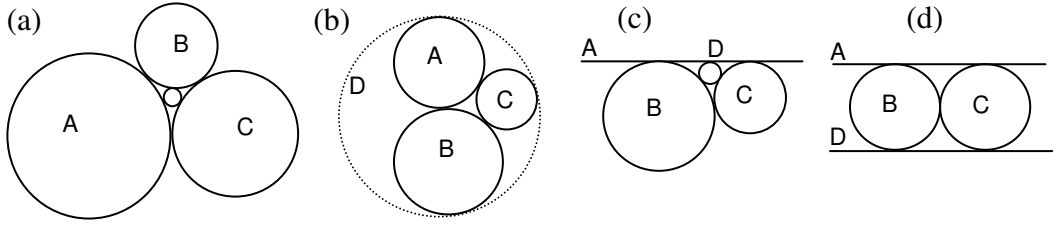

**Figure 4.1**: Four circles of the Descartes configuration.

First, consider the "all-external-tangencies" configuration (Fig. 4.1a). Construct the configuration matrix $f$ and calculate its inverse $F$:

$$f = \begin{bmatrix} -1 & 1 & 1 & 1 \\ 1 & -1 & 1 & 1 \\ 1 & 1 & -1 & 1 \\ 1 & 1 & 1 & -1 \end{bmatrix} \qquad F \equiv f^{-1} = \frac{1}{4}\begin{bmatrix} -1 & 1 & 1 & 1 \\ 1 & -1 & 1 & 1 \\ 1 & 1 & -1 & 1 \\ 1 & 1 & 1 & -1 \end{bmatrix}.$$

(Notice that $f^2 = 4\,I$; hence $f^{-1} = 1/4\,f$ readily follows). Thus our master equation $AFA^T = G$ is — after a convenient multiplication of both sides by 4 — simply:

$$\begin{bmatrix} b_1 & b_2 & b_3 & b_4 \\ \bar{b}_1 & \bar{b}_2 & \bar{b}_3 & \bar{b}_4 \\ \dot{x}_1 & \dot{x}_2 & \dot{x}_3 & \dot{x}_4 \\ \dot{y}_1 & \dot{y}_2 & \dot{y}_3 & \dot{y}_4 \end{bmatrix} \underbrace{\begin{bmatrix} -1 & 1 & 1 & 1 \\ 1 & -1 & 1 & 1 \\ 1 & 1 & -1 & 1 \\ 1 & 1 & 1 & -1 \end{bmatrix}}_{4f^{-1}} \begin{bmatrix} b_1 & \bar{b}_1 & \dot{x}_1 & \dot{y}_1 \\ b_2 & \bar{b}_2 & \dot{x}_2 & \dot{y}_2 \\ b_3 & \bar{b}_3 & \dot{x}_3 & \dot{y}_3 \\ b_4 & \bar{b}_4 & \dot{x}_4 & \dot{y}_4 \end{bmatrix} = \underbrace{\begin{bmatrix} 0 & 8 & 0 & 0 \\ 8 & 0 & 0 & 0 \\ 0 & 0 & -4 & 0 \\ 0 & 0 & 0 & -4 \end{bmatrix}}_{4g^{-1}} \qquad (4.1)$$

Matrix $D = 4f^{-1}$ (consisting of negative ones on the diagonal and ones everywhere else) is often called the Descartes matrix. The particular vector equation $\boldsymbol{b}^T F \boldsymbol{b} = G$, part of (4.1), is equivalent to the original Descartes law:

$$a^2 + b^2 + c^2 + d^2 - 2ab - 2ac - 2ad - 2bc - 2bd - 2cd = 0$$

or
$$2(a^2 + b^2 + c^2 + d^2) = (a + b + c + d)^2 \qquad (4.2)$$

In constructing matrix $f$, we have assumed that all circles are tangent externally. Let us now assume that one, say the fourth, circle contains the other three, as in Figure 4.1b. Then matrix $f$ and its inverse $F$ are:

$$f = \begin{bmatrix} -1 & 1 & 1 & -1 \\ 1 & -1 & 1 & -1 \\ 1 & 1 & -1 & -1 \\ -1 & -1 & -1 & -1 \end{bmatrix} \Rightarrow F = \frac{1}{4}\begin{bmatrix} -1 & 1 & 1 & -1 \\ 1 & -1 & 1 & -1 \\ 1 & 1 & -1 & -1 \\ -1 & -1 & -1 & -1 \end{bmatrix} = RDR,$$

where $D$ is the Descartes matrix, and $R = \text{diag}(1, 1, 1, -1)$. The Descartes formula for curvatures becomes

$$\mathbf{b}^T(RDR)\mathbf{b} = 0,$$

or

$$(R\mathbf{b})^T D\, (R\mathbf{b}) = 0.$$

This explains the convention that the circle that contains the other three circles is assumed to have a negative sign (hence the notion of "bend" as a "signed curvature"). This allows one to have a single formula for all cases in Figure 4.1.

**Corollary 4.2: (Extended Descartes Theorem, [MLW])** Four circles in Descartes configuration with possibly both types of tangencies, internal and external, satisfy the generalized Descartes theorem (4.1) if one assumes that $b$ represents bends.

(Equation (4.1) was obtained in [LMW] by different means.)

It might be a point of surprise that this convention must be replaced by other conventions, if a single formula is sought for the families of circle configurations other than the Descartes configuration. The notion of bend is **not** of a universal nature (see the following examples).

**Example 2. Beyond the Descartes configuration**

Let $a$, $b$, and $c$ be three pair-wise orthogonal circles and let $d$ be a circle tangent to each of them. There are four distinct realizations, presented here in columns:

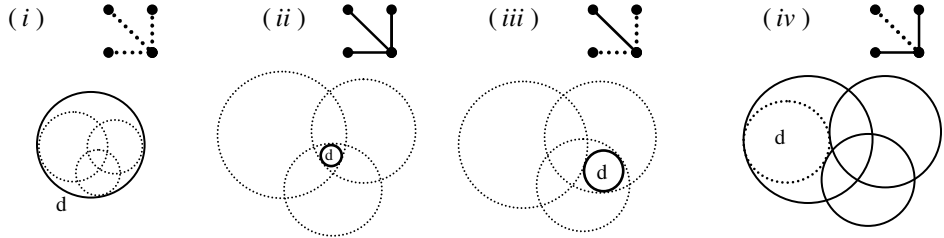

The quadratic relations may be read from these matrices. For the curvatures, one may arrive at the following concise forms:

(i)  $2(a^2 + b^2 + c^2) = (-a - b - c + d)^2$     (ii)  $2(a^2 + b^2 + c^2) = (a - b - c + d)^2$

(iii)  $2(a^2 + b^2 + c^2) = (a + b + c + d)^2$     (iv)  $2(a^2 + b^2 + c^2) = (-a + b + c + d)^2$

[Legend for the symbols: dots represent the four circles, thick lines — external tangency, dotted lines — internal tangency, and no line — orthogonality].

Since these formulas resemble each other, to account for the varying signs, the following unification may be proposed:

**Theorem 4.1:** Let three circles of curvature $a$, $b$, and $c$ be mutually orthogonal, and let a fourth circle of curvature $d$ be tangent to each of them. Then

$$2(a^2 + b^2 + c^2) = (a + b + c + d)^2$$

under the convention that if any of the mutually tangent circles, $a$, $b$ or $c$, contains or is contained in $d$, then its curvature is negative.

Note that there will always be some negative values among $a$, $b$, and $c$. If the fourth circle curvature is sought, one obtains one of the formulas $d = \pm a \pm b \pm c \pm \sqrt{2(a^2+b^2+c^2)}$.

**Example 3  Beyond the Descartes configuration again**

Let $a$ and $b$ be a pair of orthogonal circles. Add two mutually tangent circles $c$ and $d$ each of which is also tangent to $a$ and to $b$. Then we have these four cases:

(i)   (ii)   (iii)   (iv)

$$f = \begin{bmatrix} -1 & 0 & -1 & -1 \\ 0 & -1 & -1 & -1 \\ -1 & -1 & -1 & 1 \\ -1 & -1 & 1 & -1 \end{bmatrix} \quad f = \begin{bmatrix} -1 & 0 & 1 & 1 \\ 0 & -1 & -1 & -1 \\ 1 & -1 & -1 & 1 \\ 1 & -1 & 1 & -1 \end{bmatrix} \quad f = \begin{bmatrix} -1 & 0 & 1 & -1 \\ 0 & -1 & 1 & -1 \\ 1 & 1 & -1 & -1 \\ -1 & -1 & -1 & -1 \end{bmatrix} \quad f = \begin{bmatrix} -1 & 0 & 1 & 1 \\ 0 & -1 & 1 & 1 \\ 1 & 1 & -1 & 1 \\ 1 & 1 & 1 & -1 \end{bmatrix}$$

$$F = \frac{1}{8}\begin{bmatrix} -4 & 4 & -2 & -2 \\ 4 & -4 & -2 & -2 \\ -2 & -2 & -1 & 3 \\ -2 & -2 & 3 & -1 \end{bmatrix} \quad F = \frac{1}{8}\begin{bmatrix} -4 & -4 & 2 & 2 \\ -4 & -4 & -2 & -2 \\ 2 & -2 & -1 & 3 \\ 2 & -2 & 3 & -1 \end{bmatrix} \quad F = \frac{1}{8}\begin{bmatrix} -4 & 4 & 2 & -2 \\ 4 & -4 & 2 & -2 \\ 2 & 2 & -1 & -3 \\ -2 & -2 & -3 & -1 \end{bmatrix} \quad F = \frac{1}{8}\begin{bmatrix} -4 & 4 & 2 & 2 \\ 4 & -4 & 2 & 2 \\ 2 & 2 & -1 & 3 \\ 2 & 2 & 3 & -1 \end{bmatrix}$$

(i)  $2[(2a)^2 + (2b)^2 + (c-d)^2] = (-2a - 2b + c + d)^2$      (iii)  $2[(2a)^2 + (2b)^2 + (c+d)^2] = (2a + 2b + c - d)^2$

(ii) $2[(2a)^2 + (2b)^2 + (c-d)^2] = (2a - 2b + c + d)^2$       (iv)  $2[(2a)^2 + (2b)^2 + (c-d)^2] = (2a + 2b + c + d)^2$

We encounter the same problem with signs as in Example 2. (Note that the signs inside the squares may be altered when one seeks a unifying formula for all four cases). Here is a solution:

**Theorem 4.2:** Let two circles $C_1$ and $C_2$ of curvature $a$ and $b$ be mutually orthogonal, and let the second pair of circles $C_3$ and $C_4$ of curvature $c$ and $d$ be mutually tangent. Then

$$2[(2a)^2 + (2b)^2 + (c-d)^2] = (2a + 2b + c + d)^2$$

under the convention that

(i) if either of the mutually tangent circles $C_3$ and $C_4$ contains the pair
of perpendicular circles $C_1$ and $C_2$, then its curvature is negative;

(ii) if either of the perpendicular circles $C_1$ and $C_2$ contains the pair
of tangent circles $C_3$ and $C_4$, then its curvature is negative;

We now arrive at this important conclusion:

**Remark 4.3:** The case of the Descartes configuration prompted the introduction of the notion of *bend* as a signed curvature — this allowed one to write a set of formulas in a single equation (see Sec. 4). But the particular definition of "bend" — so convenient in Descartes theorem and for Apollonian gaskets — is not universal! It is native to the Descartes Theorem, and a different set of rules for the signs of "bends" may emerge in other configurations. (See also Remark 5.4).

## Example 4: A solution to Apollonius's problem

The problem of Apollonius is to find a circle that is simultaneously tangent to three geometric objects. If we choose these objects to be three circles, we may expect eight possible solutions, which differ according to the type of each of the tangencies — external versus internal (see Fig. 4.2). A geometric (constructive) solution is known (see [Coxt68]). Let us see how Theorem 3.3 provides analytic solutions.

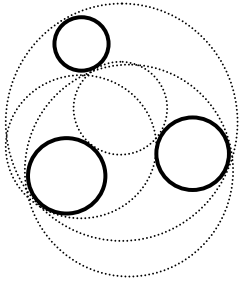

**Figure 4.2:** Four (of the eight) solutions to the Apollonian problem.

For the Apollonian problem, the configuration matrix is of the following type:

$$f = \begin{bmatrix} -1 & * & * & \pm 1 \\ * & -1 & * & \pm 1 \\ * & * & -1 & \pm 1 \\ \pm 1 & \pm 1 & \pm 1 & -1 \end{bmatrix}$$

where the last column and the last row are combinations of $+1/-1$, depending on which of the eight solutions we seek ($2^3=8$). More precisely, define coefficients that represent the configuration of the three given circles:

$$\varphi_{ij} = \frac{d_{ij}^2 - r_i^2 - r_j^2}{2 r_i r_j} = \frac{1}{2}\left(d_{ij}^2 b_i b_j - \frac{b_i^2 + b_j^2}{b_i b_j}\right),$$

where $d_{ij}$ denotes the distance between the centers of the $i$-th and $j$-th circles, $d_{ij} = d_{ji}$, $r_i$ denotes the radius of the $i$-th circle, and $b_i = 1/r_i$ is the corresponding curvature, $i,j = 1,\ldots, 3$. Then the configuration matrix is:

$$f = \begin{bmatrix} -1 & \varphi_{12} & \varphi_{13} & \varepsilon_1 \\ \varphi_{12} & -1 & \varphi_{23} & \varepsilon_2 \\ \varphi_{13} & \varphi_{23} & -1 & \varepsilon_3 \\ \varepsilon_1 & \varepsilon_2 & \varepsilon_3 & -1 \end{bmatrix}$$

where the terms $\varepsilon_n = \pm 1$ correspond to cases where the $i$-th circle *is* or *is not* contained in the circle-solution, $i = 1, 2, 3$. Finding the inverse matrix $F = f^{-1}$ and applying Theorem 3.3 gives the solution of the Apollonian problem. In particular, the quadratic equation for the curvatures is $\mathbf{b}^T F \mathbf{b} = 0$, where, as usual, $\mathbf{b} = [b_1, b_2, b_3, b_4]^T$. This needs to be solved for $b_4$, as $b_1$, $b_2$ and $b_3$ are known.

## 5. Multidimensional formulation

Consider the set $\Omega(\mathbf{E})$ of $(n-1)$-spheres in an $n$-dimensional Euclidean space $(\mathbf{E}, g_0)$, where $\mathbf{E} \cong \mathbf{R}^n$ and where $g_0$ represents the Euclidean metric. A sphere $C_r(\mathbf{p})$ of radius $r \in \mathbf{R}$ centered at $\mathbf{p} \in \mathbf{E}$ is given by equation

$$(\mathbf{x} - \mathbf{p})^2 = r^2 \tag{5.1}$$

Let $\mathbf{R}^{1,1}$ be a standard pseudo-Euclidean 2-dimensional space with metric $\omega$ given by

$$\omega = \begin{bmatrix} 0 & 1/2 \\ 1/2 & 0 \end{bmatrix}$$

Define a Minkowski space $(\mathbf{M}, g)$ by the following direct products:

$$\mathbf{M} = \mathbf{E} \oplus \mathbf{R}^{1,1} \cong \mathbf{R}^{1,n+1}, \quad g = g_0 \oplus \omega$$

**Definition 5.1:** The Pedoe map $\pi$ sends $(n-1)$-spheres into rays of the $(n+2)$-dimensional Minkowski space $\mathbf{R}^{1,n+1} \cong \mathbf{R}^{1,1} \oplus \mathbf{R}^n$ (elements of projective space):

$$\pi: \Omega(\mathbf{R}^n) \to P\mathbf{R}^{1,n+1}$$

where, for the sphere $C$ described in Eq. 8.1, the ray is spanned by vector

$$\pi(C) = \begin{cases} \text{span}\{ [ 1, r^2 - \mathbf{p}^2, \mathbf{p} ]^T \} & \text{if } C \text{ is not a plane, } r \neq 0 \\ \text{span}\{ [ 0, c/2, \mathbf{q} ]^T \} & \text{if } C \text{ is a plane } \mathbf{q} \cdot \mathbf{x} = c \\ \text{span}\{ [ 1, -\mathbf{p}^2, \mathbf{p} ]^T \} & \text{if } C \text{ is a point } \mathbf{p} \text{ (improper sphere)} \end{cases}$$

In the standard basis, the Minkowski metric $g$ in $\mathbf{M}$ and its inverse are represented by $(n+3) \times (n+3)$ matrices:

$$g = \begin{bmatrix} & 1/2 & & & \\ 1/2 & & & & \\ \hline & & -1 & & \\ & & & \ddots & \\ & & & & -1 \end{bmatrix} \qquad G = g^{-1} = \begin{bmatrix} & 2 & & & \\ 2 & & & & \\ \hline & & -1 & & \\ & & & \ddots & \\ & & & & -1 \end{bmatrix}$$

The **special Pedoe map** sends proper spheres into normed future-oriented vectors

$$\dot\pi : \Omega(\mathbf{E}) \to \mathbf{M}$$

so that (1) $|\dot\pi(C)|^2 = -1$, Vector $\mathbf{C} = \dot\pi(C)$ is called the **Pedoe vector** of sphere $S$. $\text{Im}(\dot\pi)$ lies in a hyperboloid $H \subset \mathbf{R}^{1,n+1}$ of space-like unit vectors.

Let $A = [ \dot\pi(C_1) \mid \dot\pi(C_2) \mid \ldots \mid \dot\pi(C_{n+3}) ]$ be an $(n+2) \times (n+2)$ matrix whose columns are the Pedoe vectors of a set of $n+2$ distinct $(n-1)$-spheres. The central result is:

**Theorem 5.2:** Define the **configuration matrix** for a system of $n+3$ spheres to be: $f = A^T g A$. If $f$ is invertible, then

$$AFA^T = G, \tag{5.1}$$

where $F = f^{-1}$ and $G = g^{-1}$. In particular, the curvatures and positions of $n+1$ spheres determine the remaining sphere(s) of the configuration by the quadratic relation:

$$v_i^T F v_j = G_{ij},$$

where we define $n+2$ vectors $\mathbf{v}_i$, $i = 1,\ldots,n+2$, as follows:

$\mathbf{v}_1 = \mathbf{b} = [b_1,\ldots,b_n]$, where $b_i = 1/r_i$ is the curvature of the i-th sphere;

$\mathbf{v}_2 = \bar{\mathbf{b}} = [\bar{b}_1,\ldots,\bar{b}_{n+2}]$, where $\bar{b}_i$ is the co-curvature of the i-th sphere;

$\mathbf{v}_{2+k} = \dot{\mathbf{x}}_k = [\dot{x}_{k1},\ldots,\dot{x}_{k\,n+2}]$, where $\dot{x}_{ki} = x_{ki}/r_i$ and $x_{ki}$ is the $k^{th}$ coordinate of the center of the $i^{th}$ sphere.

**Proof:** Follow the same few steps as in the proof of Theorem 3.1.

Note that the second entry of the co-curvature vector is determined by the remaining entries by the condition of normalization $\dot{\pi}(C) = -1$.

**Example 5.3:** For instance, in the case of the Descartes configuration, i.e., a system of $n+3$ mutually externally tangent $n$-spheres, the corresponding Pedoe vectors form in $\mathbf{R}^{1,\,n+1}$ a system of linearly independent vectors satisfying

$$\langle \mathbf{C}_i \cdot \mathbf{C}_j \rangle = \begin{cases} -1 & \text{if } i = j \\ 1 & \text{otherwise.} \end{cases}$$

The configuration matrix $f$ is an $(n+2) \times (n+2)$ matrix with $-1$ on the diagonal and $1$'s everywhere else, that is

$$f = N - 2I,$$

where $N$ is a matrix with $1$ in every entry, and $I$ is the unit matrix. Its inverse is $F = (1/2n)(N - nI)$; hence, after multiplying both sides of (5.1) by $2n$, we get:

$$\begin{bmatrix} b_1 & b_2 & b_3 & b_4 & \cdots \\ \bar{b}_1 & \bar{b}_2 & \bar{b}_3 & \bar{b}_4 & \cdots \\ \dot{x}_1 & \dot{x}_2 & \dot{x}_3 & \dot{x}_4 & \cdots \\ \dot{y}_1 & \dot{y}_2 & \dot{y}_3 & \dot{y}_4 & \cdots \\ \vdots & \vdots & \vdots & \vdots & \ddots \end{bmatrix} \begin{bmatrix} 1-n & 1 & 1 & 1 & \cdots \\ 1 & 1-n & 1 & 1 & \cdots \\ 1 & 1 & 1-n & 1 & \cdots \\ 1 & 1 & 1 & 1-n & \cdots \\ \vdots & \vdots & \vdots & \vdots & \ddots \end{bmatrix} \begin{bmatrix} b_1 & \bar{b}_1 & \dot{x}_1 & \dot{y}_1 & \cdots \\ b_2 & \bar{b}_2 & \dot{x}_2 & \dot{y}_2 & \cdots \\ b_3 & \bar{b}_3 & \dot{x}_3 & \dot{y}_3 & \cdots \\ b_4 & \bar{b}_4 & \dot{x}_4 & \dot{y}_4 & \cdots \\ \vdots & \vdots & \vdots & \vdots & \ddots \end{bmatrix} = 2n \begin{bmatrix} 2 & & & & \\ & 2 & & & \\ & & -1 & & \\ & & & -1 & \\ & & & & \ddots \end{bmatrix}$$

In particular, for the curvatures only, this results in

$$(b_1 + b_2 + \ldots + b_{n+2})^2 = n(b_1^2 + b_2^2 + \ldots + b_{n+2}^2),$$

which agrees with the formula discovered by Soddy [Sod] for $n = 3$, generalized to arbitrary dimension by Gossett [Gos].

**Remark 5.4 (on bends, circles and disks):** In order to account for different cases of the Descartes configuration, it is customary to introduce the notion of the *bend* of a circle: the curvature with a negative sign if the circle contains the other three circles (see e.g., [LMW]). We propose an alternative solution: replace circles by generalized disks. Any circle ($n$-sphere) is the boundary of either of two discs, one bounded and one unbounded (see Figure 5.1). The radius of the unbounded one is negative. Then each case of Descartes configuration may be represented in terms of *external* tangencies, as in Figure 5.2. The definition of the Pedoe map $\dot{\pi}$ extends to such discs: if $D$ and $D'$ are mutual complements, then we define $\dot{\pi}(D') = -\dot{\pi}(D)$.

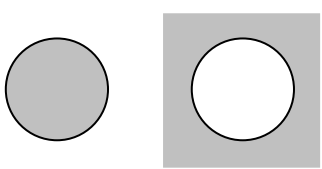

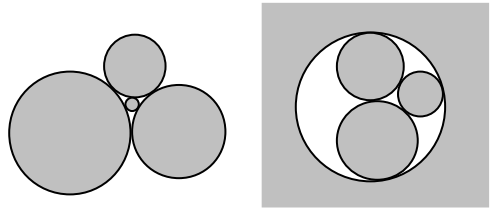

**Figure 5.1:** Bounded and unbounded disc        **Figure 5.2:** Descartes configurations reconsidered

## Acknowledgements

I would like to thanks all the participants of our Apollonian Seminar for their interest and time. This includes Philip Feinsilver, Uditha Katugampola, Xin Wang and Sankhadip Roy. Special thanks go to Alan Shoen and Greg Budzban who carefully read the manuscript and shared valuable remarks.